\documentclass[12pt]{article}
\usepackage{amsmath,amssymb,amsthm,xspace,url,picinpar,epsfig}
\usepackage{psfrag} 
\newif\ifhyper\IfFileExists{hyperref.sty}{\hypertrue}{\hyperfalse}
\hyperfalse
\hypertrue
\ifhyper\usepackage{hyperref}\fi
\def\rcs $#1: #2 ${\expandafter\def\csname rcs#1\endcsname {#2}}
\rcs $Date: 2007-11-12 $ 
\rcs $Revision: 1.10 $
\newcommand{\old}[1]{}

\newtheorem{theorem}{Theorem}[section]

\newtheorem{lemma}[theorem]{Lemma}
\newtheorem{cor}[theorem]{Corollary}
\newtheorem{prop}[theorem]{Proposition}
\newtheorem{question}[theorem]{Question}

\newcommand{\R}{{\mathbb R}}

\newcommand{\E}{{\mathbb E}}

\newcommand{\eps}{\varepsilon}

\newcommand{\sig}{\sigma}

\newcommand{\Dpgame}{{\Delta_p^\text{G}}}
\newcommand{\Dpvar}{{\Delta_p^\text{V}}}
\newcommand{\wphi}{\widetilde{\phi}}

\def\proofof#1{{ \medbreak \noindent {\bf Proof of #1.} }}
\def\dist{\text{dist}}
\def\diam{\text{diam}}

\def\BB{\overline{B}}

\begin{document}
\title{Tug of war with noise:\\ a game theoretic view of the $p$-Laplacian}
\author{Yuval Peres\thanks{Microsoft Research and U.C. Berkeley}
 \and Scott Sheffield\thanks{New York University}}
\date{} \maketitle

\begin{abstract}
Fix a bounded domain $\Omega \subset \mathbb R^d$, a continuous
function $F:\partial \Omega \rightarrow \mathbb R$, and constants
$\epsilon >0$ and $1 < p,q < \infty$ with $p^{-1} + q^{-1} = 1$. For
each $x \in \Omega$, let
$u^\epsilon(x)$ be the value for player I of the following
two-player, zero-sum game.  The initial game position is $x$.  At
each stage, a fair coin is tossed and the player who wins the toss
chooses a vector $v \in \BB(0,\epsilon)$ to add to the game
position, after which a random ``noise vector'' with
mean zero and variance $\frac{q}{p}|v|^2$ in each
orthogonal direction is also added. The game ends when the game position reaches
some $y \in\partial \Omega$, and player I's payoff is $F(y)$.

We show that (for sufficiently regular $\Omega$) as $\epsilon$ tends
to zero the functions $u^\epsilon$ converge uniformly to the unique
$p$-harmonic extension of $F$. Using a modified game (in which
$\epsilon$ gets smaller as the game position approaches $\partial
\Omega$), we prove similar statements for general bounded domains
$\Omega$ and resolutive functions $F$.

These games and their variants interpolate between the tug of war
games studied by Peres, Schramm, Sheffield, and Wilson ($p=\infty$)
and the motion-by-curvature games introduced by Spencer and studied
by Kohn and Serfaty ($p=1$). They generalize the relationship
between Brownian motion and the ordinary Laplacian and yield new
results about $p$-capacity and $p$-harmonic measure.
\end{abstract}

\section{Introduction}
Given $p>1$, a bounded domain $\Omega$ in $R^d$ and a continuous
function $F:\partial \Omega \rightarrow \mathbb R$, the
$p$-Dirichlet problem consists of finding a continuous extension $u:
\overline{\Omega} \to \mathbb R$ of $F$ which is $p$-harmonic, that
is, $u$ minimizes $\int_\Omega |\nabla u(x)|^p \,dx$ subject to the
given boundary conditions. (In general, such an extension exists
only under a regularity condition on $\Omega$, and $\int |\nabla u(x)|^p \,dx$
should be minimized over compact subsets of $\Omega$, see Proposition
\ref{viscosityweakvariationalequivalance}.) In the classical
case $p=2$, Kakutani and Doob discovered that the Dirichlet problem
can be solved by starting a Brownian motion $B$ at $x$, running it
until the hitting time $\tau$ of $\partial \Omega$, and taking
$u(x)=\E_x [F(B(\tau))]$; see \cite{MR731258} for a comprehensive
study.

In this paper, we develop an analogous interpretation of a $p$-harmonic extension $u(x)$
as the limit of the values of certain stochastic games.
The extreme case $p=\infty$ was already considered in \cite{pssw2}.

\subsection{Game definition}
Fix a bounded domain $\Omega \subset \mathbb R^d$ and a continuous
function $F:\partial \Omega \rightarrow \mathbb R$.  Let $\mu$ be a
mean zero compactly supported probability measure on $\mathbb R^d$
that is preserved by orthogonal transformations of $\mathbb R^d$
that fix the first basis vector $e_1$.  We call $\mu$ the {\em noise measure}.
(For example, $\mu$ can be the uniform distribution
on the sphere of radius $r$ in the hyperplane orthogonal to $e_1$;
 for the relation to the $p$-Laplacian described below to hold,
 $r$ should be $\sqrt{(d-1)q/p}$.)

For each $v \in \mathbb R^d$ and Borel measurable $S \subset \mathbb
R^d$, define $\mu_v(S) = \mu(\Psi^{-1}(S))$ where $\Psi$ is a
constant $c$ times some orthonormal transformation of $\mathbb R^d$,
chosen so that $\Psi(e_1)=v$.  The requirement that $\Psi(e_1)=v$
clearly implies that $c = |v|$.  In particular, if $v = c e_1$ for
some $c > 0$, then the law of a vector chosen from $\mu_v$ will be
simply the law of $c$ times a vector chosen from $\mu$. The fact
that $\mu$ is invariant under orthogonal transformations of $\mathbb
R^d$ that fix $e_1$ implies that our definition of $\mu_v$ does not
depend on the choice of $\Psi$.

Let $\alpha = 1 + \inf \{R: \mu B(0,R) = 1 \}$, where $B(z,R)$
denotes the ball of radius $R$ centered at $z$.

We now introduce a two-player zero-sum game, called {\bf tug of war
(with noise)}, played as follows. Fix an initial game state $x_0=x
\in \Omega$.  At the $k$'th turn, a fair coin is tossed, and the
player who wins the coin toss is allowed to make a move.  If
$\dist(x_{k-1}, \partial \Omega)> \alpha \epsilon$, then the moving
player chooses $v_k \in \mathbb R^d$ with $|v_k|\leq \epsilon$ and
sets $x_k = x_{k-1}+v_k + z_k$ where $z_k$ is a random ``noise
vector'' sampled from $\mu_{v_k}$. To understand the scaling,
observe that the law $\mu_{v_k}$ of the noise vector is supported on
a set of radius $(\alpha-1) |v_k|$; thus, we expect $z_k$ and $v_k$
to be of the same order of magnitude.

If $\dist(x_{k-1},
\partial \Omega) \leq \alpha \epsilon$, then the moving player
chooses an $x_k \in \partial \Omega$ with $|x_k - x_{k-1}| \leq
\alpha \epsilon$ and the game ends, with player I receiving a payoff
of $F(x_k)$ from player II.  The rules above ensure that if the game
terminates at the $k$'th step, then $x_j \in \Omega$ for all $0 \leq
j < k$ and $x_k \in \partial \Omega$. Both players receive a payoff
of zero if the game never terminates.
\begin{figure}[htbp]\label{fig1}
\begin{center}
\psfrag{xk?}{${x_k}?$} \psfrag{zk?}{${z_k}?$}
\psfrag{xk-1}{${x_{k-1}}$} \psfrag{vk}{${v_{k}}$}
\psfrag{eps}{$\epsilon$} \psfrag{reps}{$r\epsilon$}
\epsfig{file=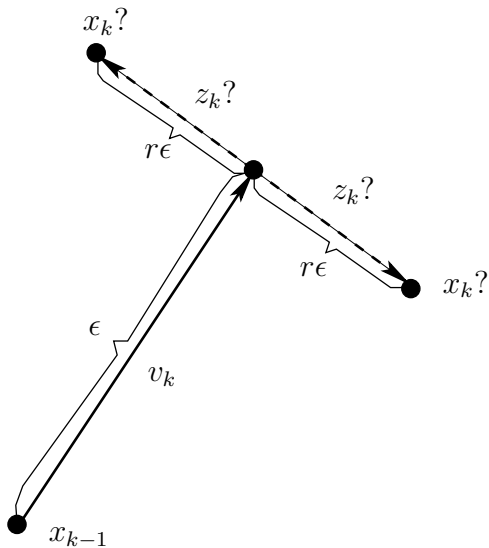, width=.40 \textwidth} \caption{A move of tug
of war with noise in dimension 2, for the noise distribution $\mu$
given by $\mu\{(0,r)\}=\mu\{(0,-r)\}=1/2$. The player who wins the
coin toss adds a vector $v_k$ of length at most $\epsilon$ to the
the game position $x_{k-1}$, and then a random noise vector $z_k$
with law $\mu_{v_k}$ (and magnitude $r|v_k|$) is added to produce $x_k$.
In the figure $|v_k|=\epsilon$.}
\end{center}
\end{figure}
Tug of war in a domain of $\mathbb R^d$ without noise (i.e., with
$\mu$ supported at the origin) was introduced in \cite{pssw2}, where
it was used to give uniqueness results for PDEs involving the
infinity Laplacian and also extended to solve the optimal Lipschitz
extension problem on general length spaces. Other games with random
turn order were introduced earlier \cite{pssw, MR1427981}.

\subsection{Main results}
Fix the game parameters $F$, $\Omega$, and $\mu$ above, and denote by
$u^\epsilon_1(x)$ the supremum over all player I (measurable, pure)
{\bf strategies} (defined precisely in Section
\ref{strategydefinitionsection}) of the infimum over all player II
strategies of the expected payoff for player I when both players
adopt these strategies. Denote by $u^\epsilon_2(x)$ the infimum over
all player II strategies of the supremum over all player I
strategies of the expected payoff for player I when both players
adopt these strategies. It is clear that   $u^\epsilon_1 \le u^\epsilon_2$.  When $u^\epsilon_1 = u^\epsilon_2$ we say
that the game {\bf has value} $u^\epsilon = u^\epsilon_1 =
u^\epsilon_2$.

A function $u$ is {\bf $p$-harmonic} in a domain $\Omega \subset
\mathbb R^d$ if $u$ is continuous on $\Omega$ and for every bounded
subdomain $\Omega_0$ of $\Omega$ such that $\overline \Omega_0
\subset \Omega$, the {\bf $p$-energy} $\int_{\Omega_0} |\nabla
u(x)|^pdx$ of $u$ in $\Omega_0$ is finite and as small as possible,
given the values of $u$ on $\partial \Omega_0$. (We give equivalent
viscosity and weak $p$-Laplacian-based definitions in Section
\ref{analysisdefinitionsection}, where we will also define the terms
``$p$-superharmonic'' and ``$p$-subharmonic.'') Note that the continuity assumption can be replaced by a weaker regularity
assumption, and that
$p$-harmonic functions are the same as harmonic functions when
$p=2$.

Let $\Pi_i$ denote projection to the $i$th coordinate.
The covariance matrix $$C=\Bigl(\int \Pi_i(x) \Pi_j(x) \, d\mu(x) \Bigr)_{i,j=1}^d$$
 of the noise measure $\mu$ is necessarily diagonal with
$C_{i,i} = C_{j,j}$ for all $2 \leq i,j \leq d$.
The main result of this paper is that as $\epsilon$ tends to zero,
the functions $u^{\epsilon}_1$ and $u^\epsilon_2$ converge uniformly
to the unique $p$-harmonic extension of $F$ (at least when $F$ is
continuous and $\Omega$ is sufficiently regular), where $p$ is
determined as follows: Given $\mu$, we
define a constant $p = p(\mu) \in [1, \infty]$ by $p =
\frac{C_{1,1}+C_{2,2}+1}{C_{2,2}}$. Equivalently, $p$ is such that
for some $\beta
> 0$, we have $C_{1,1}+1 = \beta
q^{-1}$ and $C_{i,i} = \beta p^{-1}$ for $2 \le i \le d$, where $p^{-1} + q^{-1} = 1$.
This $p$ is chosen so that if one player always chooses the vector
$v_k$ to be some $v$ with $|v|=\epsilon$, and the other player
always chooses $v_k$ to be $-v$, then for each $k$ the variance of
$x_k - x_{k-1}$ is proportional to $q^{-1}\epsilon^2$ in the
direction of $v$ and $p^{-1}\epsilon^2$ in each direction orthogonal
to $v$.  In other words, $q^{-1}$ and $p^{-1}$ are the relative
sizes of the parallel and perpendicular variance when the players
``tug'' in opposite directions.  The case $p=\infty$ corresponds to
purely parallel variance and $p=1$ to purely perpendicular variance.

We now describe the required regularity condition on $\Omega$. For a
given noise measure $\mu$ and $p=p(\mu)$, a point $y \in
\partial \Omega$ is called a {\em game-regular boundary point} of
$\Omega$ if whenever the game starts near that $y$, player I has a
strategy for making the game terminate near $y$ with high
probability. More precisely, $y$ is {\bf game-regular} if for every
$\delta > 0$ and $\eta
> 0$ there exists a $\delta_0$ and $\epsilon_0$ such that for every
$x_0$ with $|x_0 - y| < \delta_0$ and $\epsilon < \epsilon_0$,
player I has a strategy that guarantees that an $\epsilon$-step game
started at $x_0$ will terminate at a point on $\partial \Omega \cap
B(y,\delta)$ with probability at least $1 - \eta$. We say $\Omega$
itself is {\bf game-regular} if every $y \in
\partial \Omega$ is game-regular. We say that a point $y \in \partial \Omega$ satisfies the {\bf
cone property} if there is a neighborhood $U$ of $y$ and a cone
$\mathcal C$ with tip $y$ such that $\mathcal C \cap U$ is disjoint
from $\Omega$.  We say that $\Omega$ satisfies the cone property if every
$y \in \partial \Omega$ satisfies the cone property.  The following
results will be proved in Section \ref{proofsection}.

\begin{prop} \label{coneimpliesgameregular} Let $\Omega$ be a bounded domain in $\mathbb R^d$
and fix a noise measure $\mu$ and $p=p(\mu)\in(1,\infty)$ as above. Then
\begin{description}
\item{{\bf (i)}} If $p
> d$ then $\Omega$ is game-regular.
\item{{\bf (ii)}} Every $y
\in
\partial \Omega$ that satisfies the cone property is game-regular.
\item{{\bf (iii)}} If $d=2$ and $\Omega$ is simply connected then $\Omega$ is
game-regular.
\end{description}
\end{prop}

\begin{theorem} \label{regularandcontinuousconvergence}
Let $\Omega$ be a bounded domain in $\mathbb R^d$.
Fix a noise measure $\mu$ and $p=p(\mu)\in(1,\infty)$.
\begin{description}
\item{{\bf (i)}}
Suppose $\Omega$ is game-regular and $F:
\partial \Omega \rightarrow \mathbb R$ is continuous.  Then as $\epsilon \rightarrow 0$,
the game values $u^\epsilon_1$ and $u^\epsilon_2$ converge uniformly
to the unique $p$-harmonic function $u$ that extends continuously to
$F$ on $\partial \Omega$.
\item{{\bf (ii)}} Conversely, if $\Omega$ is not game
regular, then there exists a continuous function $F:\partial \Omega \rightarrow \mathbb R$
 for which $u^\epsilon_1$ does not
converge uniformly to a function that extends continuously to
$\partial \Omega$.
\end{description}
\end{theorem}

Section \ref{variantsection} extends the above results to two
natural variants of tug of war.  The first is the same as the game
described above except that the players alternate turns instead of
deciding turn order with coin tosses.  Alternating turn tug of war
without noise is not interesting (since either player may choose a
strategy of always undoing the other player's moves, thereby
preventing the game from terminating), but the values in alternating
turn tug of war with noise converge to $p$-harmonic functions for an
appropriate choice of $p$.  The second variant interpolates between
tug of war without noise ($p = \infty$) and certain optimal control
processes introduced by Spencer nearly thirty years ago
\cite{MR0526057} and studied in detail by Kohn and Serfaty
\cite{MR2200259} ($p=1$).

Analogous results for discontinuous boundary conditions and
non-regular domains are presented in Section
\ref{discontinuousboundarysection}.
Nonmeasurable strategies are briefly discussed in
Section \ref{valueexistencesection},
and some unsolved problems are collected in the final section.

\subsection{Motivation}

In the case $p= \infty$, the tug of war games without noise were
used to prove new uniqueness results for
$\Delta_\infty u = g$, to solve the optimal Lipschitz extension
problem on general length spaces, and to give various bounds on
infinity harmonic measure \cite{pssw2}.  We hope that the results of
this paper will be similarly useful to the study of $p$-harmonic
measure and $p$-capacity.   In Section \ref{harmonicmeasuresection}
we will see that some of the bounds in \cite{pssw2} translate to the
$p$-harmonic case as well.  In particular, the question of how the
$p$-harmonic measure of a $\delta$ neighborhood of a Cantor set
decays as $\delta$ tends to zero can be easily addressed using game
theoretic arguments.

Furthermore, although the games in \cite{pssw} were applied to solve
problems in analysis, the authors of \cite{pssw} (including the
current authors) were originally motivated by the games themselves.
We expect that tug of war games will have applications to political
and economical modeling. They are natural models for situations in
which the relevant state space is well summarized by finitely many
parameters and in which opposing parties continually seek to improve
their positions though incremental ``tugs.''  They are related to
many of the differential games (in which players choose drift and
diffusion terms at each point in a domain) used in economic modeling,
but they are particularly simple in that the move sets are player-symmetric, and
each player's allowed set of incremental moves is independent of
anything the other player does. Barron, Evans, and Jensen have
recently developed some continuous time variants of tug of war without noise in
which both players specify a {\em control flow} on the entire domain
$\Omega$ in advance, and then the game position alternates randomly
between evolving according to one player's control and evolving
according to the other player's control until it reaches the
boundary \cite{BEJ}.

\subsection{More definitions and background about $p$-harmonicity}
\label{analysisdefinitionsection}

On the space of real-valued functions on a subset $\Omega \subset
\mathbb R^d$, we define the {\bf infinity Laplacian} operator
$\Delta_\infty$ by $\Delta_\infty u = |\nabla u|^{-2} \sum_{i,j} u_i
u_j u_{ij}$. (Here $u_i = \frac{\partial u}{\partial x_i}$.)

When $u$ is twice differentiable and $\nabla u \not = 0$, the
infinity Laplacian is the second derivative of $u$ in the gradient
direction. For more on the infinity Laplacian, see the survey paper
\cite{MR2083637}. We similarly define the {\bf $1$-Laplacian}
$\Delta_1 = \Delta - \Delta_\infty$ where $\Delta u = \sum_{i=1}^d
u_{ii}$ is the ordinary Laplacian. In other words, $\Delta_1 u$ is
the sum of the second derivatives in each of the $d-1$ directions
orthogonal to $\nabla u$.

For $1 < p < \infty$, a function that minimizes $\int |\nabla
u(x)|^pdx$ (given its boundary values) solves the {\bf Euler
Lagrange equation}
\begin{equation} \label{eulerlagrange}\text{div}(|\nabla
u|^{p-2}\nabla u) = 0.\end{equation}  When $u$ is smooth, we can
write the left hand side as \begin{eqnarray*} \Dpvar u & := &
 \sum_{i=1}^d
\frac{\partial}{\partial x_i} \left(|\nabla u|^{p-2} u_i \right) \\
& = & \sum_{i=1}^d |\nabla u|^{p-2} u_{ii} + \sum_{i=1}^d u_i
\frac{\partial}{\partial x_i} \left((|\nabla u|^2)^{(p-2)/2}\right) \\
& = & |\nabla u|^{p-2} \Delta u + \sum_{i=1}^d u_i \frac{p-2}{2}
|\nabla u|^{p-4}\frac{\partial}{\partial x_i} |\nabla u|^2 \\
& = & |\nabla u|^{p-2} \Delta u + (p-2)|\nabla u|^{p-4} \sum_{i=1}^d
\sum_{j=1}^d u_i u_{ij} u_j \\
& = & p |\nabla u|^{p-2} \left(p^{-1} \Delta_1 +
q^{-1}\Delta_\infty\right)u.
\end{eqnarray*} If $\nabla u \not = 0$, then the last expression vanishes precisely
when $p^{-1} \Delta_1 + q^{-1}\Delta_\infty = 0$.  Thus, for general
$1 \leq p,q \leq \infty$ (with $p^{-1} + q^{-1} = 1$), it is natural
to define an operator $\Dpgame$ to be the convex combination of
$\Delta_1$ and $\Delta_\infty$ given by $\Dpgame := p^{-1} \Delta_1
+ q^{-1}\Delta_\infty$.  Since $\Dpgame$ and $\Dpvar$ differ only by
a normalization factor of $p |\nabla u|^{p-2}$, the equation
$\Dpgame u = 0$ is equivalent in the classical sense to the
Euler-Lagrange equation $\Dpvar u = 0$ (provided $\nabla u \not =
0$).  It is also equivalent in the weak sense as well as the
viscosity sense, as discussed below. Note that if $\nabla u(x) = 0$,
then $\Dpgame u(x)$ is undefined in the classical sense (as is
$\Dpvar u (x)$ when $p<2$).

In much of the literature, the $p$-Laplacian is set to zero, so the
distinction between $\Dpvar$ and $\Dpgame$ is irrelevant. The
equations $\Dpvar u = g$ for non-zero constant functions $g$ have a
natural variational interpretation in terms of minimizing the
$p$-energy of $u$ conditioned on the volume $\int_\Omega u(x) dx$
bounded underneath $u$. However, the solution $u$ of $\Dpgame u =
g=-2f/\beta$ has a natural game-theoretic interpretation as the
limiting value of a modified tug of war game in which in addition to
the payoff player I receives when the game terminates, player I
receives a ``running payoff'' of size $f(x_k) \epsilon^2$ at the
$k$th step of the game. (See \cite{pssw2} for complete details in
the case $p = \infty$.) The normalization $\Dpgame$ also has the
aesthetic advantage of being a convex combination of $\Delta_1$ and
$\Delta_\infty$. Throughout the remainder of this paper, when we use
the term {\bf $p$-Laplacian}, or write $\Delta_p$, we will always
mean $\Dpgame$.


We say that $u$ is a {\bf viscosity subsolution} to $\Delta_p (\cdot) =
g$ in an open set $U$ if $u$ is upper semi-continuous and for every
$x \in U$ and $C^2$ function $\phi$ on a neighborhood of $x$ such
that
\begin{enumerate}
\item $\phi(x) = u(x)$,
\item $u \leq \phi$ in a neighborhood of $x$,
\item $\nabla \phi(x) \not = 0$,
\end{enumerate} we have $\Delta_p \phi(x) \geq g(x)$.  We say $u$ is a {\bf viscosity
supersolution} to $\Delta_p (\cdot) = g$ if $-u$ is a viscosity subsolution to $\Delta_p (\cdot) =-g$.  We make similar
definitions for $\Dpvar$ and observe that $u$ is a viscosity
subsolution to $\Delta_p u= 0$ if and only if it is a viscosity
subsolution to $\Dpvar u = 0$.  (In the latter case, when $p \geq
2$, the relevant comparison functions $\phi$ are those that have
$\Dpvar \phi(x) > 0$ at $x$, which in turn implies $\nabla \phi \not
= 0$ and $\Delta_p \phi(x)> 0$.)

Both \cite{MR1871417} and \cite{MR1207810} define $p$-harmonic
functions to be weak solutions to $\Dpvar u = 0$: precisely, a
function $u$ on a domain $\Omega$ is called {\bf $p$-harmonic} in
$\Omega$ (for $1 < p < \infty$) if it is continuous and belongs to
the Sobolev space $W^{1,p}_{\text{loc}}(\Omega)$ (i.e., $\nabla u
\in L^p_{\text{loc}}(\Omega)$) and

\begin{equation} \label{euler} \int_{\Omega} |\nabla u|^{p-2} ( \nabla u, \nabla \phi ) dx =
0 \end{equation} for every $\phi \in C^{\infty}_0(\Omega)$.  (This
is the distributional form of (\ref{eulerlagrange}).) Say that
$u:\Omega \rightarrow (-\infty,\infty]$ is $p$-superhamonic if $u$
is lower semicontinuous, $u \not = \infty$, and $u$ satisfies the
following comparison principle on each bounded subdomain $\Omega_0$
with $\overline \Omega_0 \subset \Omega$: if $h \in
C(\overline{\Omega}_0)$ is $p$-harmonic in $\Omega_0$ and $u \geq h$
on $\partial \Omega_0$, then $u \geq h$ in $\Omega_0$. We say $u$ is
$p$-subharmonic if $-u$ is $p$-superharmonic.

We now cite the following:

\begin{prop} \label{viscosityweakvariationalequivalance} When $1 < p < \infty$, a function is
$p$-superharmonic in $\Omega$ if and only if it is a viscosity
$p$-supersolution.  Moreover, $u$ is $p$-harmonic (equivalently, a
viscosity $p$-solution) if and only if $u$ is a continuous solution
to the variational problem, i.e., for every subdomain $\Omega_0$
with compact closure $\overline \Omega_0 \subset \Omega$, the
integral $\int_{\Omega_0} |\nabla u|^p$ is finite and as small as
possible, given the values of $u$ on $\partial \Omega_0$.
\end{prop}

The viscosity and weak equivalence was proved in \cite{MR1871417}.
The equivalence of (\ref{euler}) with the variational problem is
classical Euler-Lagrange theory; see, e.g., the reference text
\cite{MR1207810} on non-linear potential theory. Recall that
$W^{1,p}_0(\Omega)$ is the closure in $W^{1,p}$ of the smooth
functions supported on compact subsets of $\Omega$.  From Chapter 3
and Chapter 5 of \cite{MR1207810} we cite the following:

\begin{prop} If
$\Omega \subset \mathbb R^d$ is bounded and $\theta \in
W^{1,p}(\Omega)$, then there exists a unique $p$-harmonic function
$u$ in $\Omega$ such that $u - \theta \in W^{1,p}_0(\Omega)$;
moreover, this $u$ minimizes the integral $\int_\Omega |\nabla u|^p$
among all functions $v \in W^{1,p}(\Omega)$ for which $v-\theta \in
W^{1,p}_0(\Omega)$.
\end{prop}

\begin{prop} \label{pharmonicunique} If $\Omega$ is bounded and $\phi \in C(\overline \Omega)$ is $p$-harmonic
in $\Omega$, then any $p$-harmonic function which extends
continuously to the boundary of $\Omega$ and agrees with $\phi$ on
the boundary of $\Omega$ is equal to $\phi$.
\end{prop}

We also quote the following smoothness result, which was proved by DiBendetto~\cite{ben} and independently by Tolksdorf~\cite{Tolk}, who extended earlier work by Morrey, Uhlenbeck, Evans and others.  
\begin{prop} \label{pharmonicsmoothness}
If $u$ is $p$-harmonic in $\Omega$, then it is everywhere
differentiable in $\Omega$ and real analytic wherever $\nabla u \not
= 0$.  Moreover, $u$ has a H\"{o}lder continuous gradient---i.e., $u
\in C^{1,\gamma}(\Omega)$ for some $\gamma > 0$.
\end{prop}

\section{Proofs} \label{proofsection}
\subsection{Strategy definition} \label{strategydefinitionsection}

When two players play tug of war, we define the {\bf history up to
step $k$} of the game to be the sequence $h_k= \{x_0,v_1, x_1,v_2,
x_2, \ldots, v_k, x_k\}$.   (If the game terminated at time $j < k$,
we set $v_m = 0$ and $x_m = x_j$ for $m \geq j$.)  This $h_k$
belongs to the space $H_k = \Omega \times (\BB(0,\epsilon)
\times \overline \Omega)^k$.  The complete {\bf history space}
$H_\infty$ is the set of all infinite game position sequences $h=
\{x_0, v_1, x_1, v_2, x_2, \ldots \}$.

We endow $H_\infty$ with the product topology. It is easy to see
that if $F:\partial \Omega \to \mathbb R$ is Borel measurable, then
the payoff for player I is a Borel measurable function on
$H_\infty$.  A {\bf (measurable pure) strategy} is a sequence of
Borel measurable maps from $H_k$ to $\overline B(0,\epsilon)$,
giving the move a player would make at the $k$th step of the game as
a function of the game history.
 A pair of strategies $\sigma=(S_I, S_{II})$ (where $S_I$ is a strategy for player I and $S_{II}$ is a strategy for player II) and a starting point $x$ determine a unique probability measure $\mathbf{P}_x$ in H1 (constructed using Kolmogorov's
extension theorem). Denoting the corresponding expectation by $\mathbf{E}_x$,
the expected payoff for this strategy pair is
$ V (S_I, S_{II}) = \mathbf{E}_x  \Bigl[F(x_\tau) \mathbf{1}_{\tau<\infty}\Bigr]$
where $\tau$ is the exit time of $\Omega$.
We then have
$u_1^\eps = \sup_{S_I} \inf_{S_{II}}  V (S_I, S_{II}) $
and
$u_2^\eps =  \inf_{S_{II}} \sup_{S_I} V (S_I, S_{II}) $.

The restriction to Borel measurable strategies is necessary since
$\mu$ is a Borel measure.  It is also natural, since measurable
strategies are in some sense the only ones that can be implemented
in a constructive way.  (We mention some relaxations of the
measurability condition in Section \ref{valueexistencesection}.)

\subsection{Smooth case, non-vanishing gradient}
If $x$
is a matrix or column vector, denote by $x^T$ the transpose of $x$.
Suppose that $z$ has law $\mu$ and $z_1=\Psi(z)$ has law $\mu_v$
(where $\Psi$ is a rotation multiplied by $|v|$ and $\Psi(e_1)=v)$.
As noted above, $z$ has mean zero and covariance matrix
$$C^{z} = (\beta q^{-1}-1) e_1 e_1^T + \beta p^{-1} (I - e_1 e_1^T).$$
 Therefore, $z_1$ has mean zero and covariance matrix
$$C^{z_1} = (\beta q^{-1}-1) v v^T + \beta p^{-1} (|v|^2I - vv^T).$$

Next, we need some basic lemmas about quadratic functions.
Let $\phi(x) = x^TAx + (\xi,x)$ (for $\xi \in \R^d \setminus\{0\}$) be a quadratic
function. Without loss of generality, $A$ is a symmetric matrix, since only the quadratic form induced
by A is used. Suppose that $x_0 = 0$ and let $\psi(v)$ be the expected
value of $\phi(x_1) = \phi(v_1 + z_1)$ if player I chooses $v_1 =
v$.
Let $z_1^i$ denote the $i$th component of $z_1$. Since $\mathbb E
z_1 = 0$, we have $$\psi(v) := \mathbb E \phi(v+z_1) = \phi(v) +
\sum_{i,j} A_{i,j} \mathbb E z_1^i z_1^j = \phi(v)+\sum_{i,j} A_{i,j}
C^{z_1}_{i,j}.$$ Note that $\sum_{i,j} (v v^T)_{i,j} A_{i,j} = v^T A
v$ and $\sum_{i,j} I_{i,j} A_{i,j} = \text{Tr} A$. Therefore

\begin{eqnarray*}
\psi(v) &= &\phi(v) + (\beta q^{-1}-1) v^T A v + \beta p^{-1}
(|v|^2\text{Tr}
A - v^T A v) \\
&= &\phi(v) + (\beta q^{-1} - \beta p^{-1}-1)v^TAv + \beta p^{-1} |v|^2\text{Tr}A\\
& =& (\beta q^{-1} - \beta p^{-1}) v^TAv + \beta p^{-1}
|v|^2\text{Tr} A + (\xi, v),
\end{eqnarray*} whence $\psi(v) = (\xi, v) + v^T B v$ where $B =
(\beta q^{-1} - \beta p^{-1})A + \beta p^{-1}(\text{Tr}A) I$. In the
following, we use the $L^2$ operator norm on matrices given by $\|B\| =
\sup_{|v|=1} |B v|$.

\begin{lemma} \label{locationofmaximizer}


Fix $\xi \in \R^d \setminus\{0\}$, and let $\zeta = \frac{4\|B\|}{|\xi|}.$  If $\epsilon < \zeta^{-1}$,
then any $v \in \overline{B}(0,\epsilon)$ for which $\psi(v) = (\xi,
v) + v^T B v$ is maximal within $\overline{B}(0, \epsilon)$,
satisfies $|v|= \epsilon$ and $\vline\, v - (\epsilon |\xi|^{-1}\xi)
\vline\, \leq \zeta \epsilon^2$.
\end{lemma}

\begin{proof}
We first observe that the maximum of $\psi$ within
$\BB(0,\epsilon)$ is obtained on the boundary $\partial
B(0,\epsilon)$: Otherwise $\nabla \psi(v) = \xi + 2Bv$ must vanish
at the maximum, but $$|2Bv| \leq 2\|B\|\epsilon < 2\|B\|\zeta^{-1} <
|\xi|.
$$

Let $w = \epsilon|\xi|^{-1}\xi$.  Then $\psi(v) \geq \psi(w) =
\epsilon|\xi|+w^TBw$.  This means
$$(\xi, v) \geq \epsilon|\xi| + w^TBw - v^TBv \geq \epsilon |\xi| -
2\|B\|\epsilon|v-w|.$$ Multiplying by $\epsilon|\xi|^{-1}$ gives
\begin{equation} \label{scalarequation}(w,v) \geq \epsilon^2 - \epsilon^2\frac{\zeta}{2}|v-w|,\end{equation}
whence
$$|v-w|^2 = 2\epsilon^2 - 2(w,v) \leq \zeta |v-w| \epsilon^2.$$

\end{proof}

\begin{lemma} \label{maximalgainperturn}
Given $v$, $\psi$, $B$, and $\xi$ as in the statement Lemma
\ref{locationofmaximizer}, we have the following approximation of
$\psi(v)$:
$$\vline \,\psi (v) - \epsilon |\xi| - \Delta_\infty\psi (0) \epsilon^2/2 \,\vline \leq
16 \frac{\|B\|^2}{|\xi|} \epsilon^3,$$ for all $\epsilon < 1$.
\end{lemma}

\begin{proof}
Since $\Delta_\infty \psi(0)$ is the second derivative of $\psi$ in
the $\xi$ direction, the vector $w = \epsilon |\xi|^{-1}\xi $ satisfies
$\psi(w) = \epsilon |\xi| + \Delta_\infty \psi(0) \epsilon^2/2.$ Let
$\zeta$ be as in Lemma \ref{locationofmaximizer}. Recall that $|v| =
|w| = \epsilon$ and $|w - v| < \zeta \epsilon^2$. By
(\ref{scalarequation}), $$(w,w-v) \leq \frac{\epsilon^2 \zeta}{2}
|v-w| \leq \frac{\epsilon^4 \zeta^2}{2}.$$  Therefore,
\begin{eqnarray*}
|\psi(v)-\psi(w)| &\leq& \frac{|\xi|}{\epsilon} |(w,w-v)|+
2\|B\|\epsilon|v-w| \\
&\leq& \frac{|\xi|\epsilon^3\zeta^2}{2} + 2\|B\| \zeta \epsilon^3 =
\frac{16\|B\|^2}{|\xi|}\epsilon^3.
\end{eqnarray*}

\end{proof}

If $\|B\|$ is bounded above and $|\xi|$ is bounded below, then Lemma
\ref{maximalgainperturn} implies that the expected change in $\phi$
when player I makes an optimal move is the same as the expected
change when player I moves in the gradient direction, up to an error
of $O(\epsilon^3)$. These estimates deteriorate badly when $|\xi|$
tends to zero and $\|B\|$ gets large.  However, they are enough to
prove the following:

\begin{lemma} \label{onestepexpectationbound}
Let $\phi(x)=x^TAx+(\xi,x)$.  Fix $k \ge 0$ and suppose that player
I's strategy in move $k+1$ is to tug distance $\epsilon$ in the
gradient direction of $\phi$ if player I wins the coin toss. Then
regardless of player II's strategy, we have
$$\mathbb E [\phi(x_{k+1})|h_k] \ge \phi(x_k) + \frac{\beta}{2}
\Delta_p \phi(x_k)\epsilon^2- M\epsilon^3,$$ where $$h_k = \{x_0,
v_1, x_1, v_2, x_2, \ldots, v_k, x_k\}$$ is the game history up to
step $k$ and $M = \frac{16 \beta (d+1)\|A\|^2}{|\xi|}$.
\end{lemma}

\begin{proof} Recall that $\psi(v) = (\xi, v) + v^T B v$ where $B =
(\beta q^{-1} - \beta p^{-1})A + \beta p^{-1}(\text{Tr}A) I$. Thus
$$\Delta_\infty \psi(0) = (\beta q^{-1} - \beta p^{-1})\Delta_\infty
\phi + \beta p^{-1} \Delta \phi = \beta \Delta_p \phi(0).$$ Lemma
\ref{maximalgainperturn} implies that if player I adopts the
strategy in the statement, then
$$\mathbb E [\phi(x_{k+1})|h_k^I] \ge \phi(x_k) + \epsilon|\xi| + \frac{\beta}{2}
\Delta_p \phi(x_k)\epsilon^2- M\epsilon^3,$$ where $h_k^I$ indicates
that in move $k+1$, player I won the coin toss and chooses $v_{k+1}$, whence the conditional expectation is simply integrating with respect to $\mu_{v_{k+1}}$. Minimizing $\psi$ over $\overline{B}(0, \epsilon)$
is the same as maximizing $-\psi$ there. Now $-\psi$ is obtained from $\psi$ by replacing $\xi$ and $B$ by $-\xi$ and $-B$, respectively. Therefore, Lemma
\ref{maximalgainperturn} also implies that regardless of player II's
strategy, we have $$\mathbb E [\phi(x_{k+1}) |h_k^{II}] \geq
\phi(x_k) -\epsilon|\xi| + \frac{\beta}{2} \Delta_p
\phi(x_k)\epsilon^2 -M\epsilon^3 ,$$
where the conditional expectation is integrating with respect to $\mu_{v_{k+1}}$ as before, but now
$v_{k+1}$ was chosen by player $II$.
Averaging the last two displayed
equations proves the lemma.
\end{proof}

Given a differentiable function $u$ with non-vanishing gradient, the
{\bf gradient strategy} for player I is to always take $$v_k =
\epsilon |\nabla u(x_{k-1})|^{-1} \nabla u(x_{k-1})$$ at every step
of the game regardless of what player II does.  The {\bf gradient
strategy} for player II is to always take $v_k =-\epsilon |\nabla
u(x_{k-1})|^{-1} \nabla u(x_{k-1})$. We can interpret Lemma
\ref{onestepexpectationbound} as a bound on the expected change in
$u$ when the gradient strategy is employed.

\begin{theorem} \label{pharmoniclimitsmoothcase}
Suppose that $u$ is $p$-harmonic on a domain $\Omega'$ and $\nabla u
\ne 0$ throughout $\Omega'$. Let $\Omega$ be a domain whose closure
is a compact subset of $\Omega'$. Fix a noise measure $\mu$. For $x
\in \Omega$, let $u^\epsilon_1(x)$ be the value for player I of the
game in $\Omega$, started at $x_0=x$, with boundary conditions $F =
u$ on $\partial \Omega$. Then the functions $u^\epsilon_1$ converge
uniformly to $u$. In fact, $|u^\epsilon_1 - u|_\infty =
O(\epsilon)$. (The implied constant may depend on $\mu$.)
\end{theorem}

\begin{proof}
We may assume (adding a constant if necessary) that $u > 0$ on
$\overline{\Omega}$.  By Proposition \ref{pharmonicsmoothness}, $u$
is real analytic in $\Omega'$, and in particular its derivatives
restricted to $\Omega$ are all bounded.  We thus have the Taylor
expansion at any $y \in \Omega$ given by $u(x) = u(y) + (\nabla u,
(x-y)) + \frac{1}{2} (x-y)^T (D^2 u) (x-y) + O((x-y)^3)$. Thus, at
every $x$, the function $\phi(x) = (\nabla u, y) + (x-y)^T D^2 u
(x-y)$ approximates $u$ up to an error of $O(\epsilon^3)$ within the
ball of radius $2(\alpha + 1)\epsilon$ of $y$.  The bounds in Lemma
\ref{onestepexpectationbound} then imply that player I, by adopting
the gradient strategy, can ensure that for some constant $c$ the
sequence $M_k = u(x_k) +  ck \epsilon^3$ is a submartingale.

The non-vanishing gradient and compactness imply that for some $c_1
> 0$, we have $u(x_k) - u(x_{k-1}) \geq c_1 \epsilon$ provided
player I wins the coin toss and $\epsilon$ is small enough.
Therefore $\mathbb E[(M_k - M_{k-1})^2 | h_{k-1}] \ge \frac{1}{2} c_1^2 \epsilon^2$.
 Consequently, the
difference $M_k^2 - M_0^2 - c_2 \epsilon^2 k$, where $c_2 =
c_1^2/2$, is a submartingale, because
$$\mathbb E [M_k^2 - M_{k-1}^2 | h_{k-1}] = \mathbb E[(M_k - M_{k-1})^2 | h_{k-1}] + \mathbb E[2(M_k
- M_{k-1})M_{k-1}| h_{k-1}],$$ and the second term on the right hand
side is non-negative.

\medskip
Suppose that  player I adopts the gradient strategy and player II adopts some
arbitrary strtaegy $S_{II}$.
Use $\mathbb P_x$ to denote probability when the initial game
position is $x_0 =x$. Then
have
\begin{eqnarray*} c_2k
\mathbb P_x[ \tau \geq k] &\leq & c_2\mathbb E_x[\tau \wedge k] \\
&\leq& \epsilon^{-2} \mathbb E_x[M^2_{\tau \wedge k}
- M_0^2]\\
& \leq &\epsilon^{-2}(8\|u\|_\infty^2 + 2c^2k^2\epsilon^6 ),\\
\end{eqnarray*}
so $\mathbb P_x(\tau > a\epsilon^{-2}) \leq 1/2$ for some $a$
(independent of $\epsilon$).

Since this holds uniformly in $\epsilon$ and $x$, iteration shows
that $\frac{\tau}{a \epsilon^{-2}}$ can be bounded by a
$\text{geometric}(1/2)$ random variable.  Thus, $\mathbb E_x[\tau] =
O(\epsilon^{-2})$.  Applying the optional stopping time theorem for
submartingales, $$u(x) = \mathbb E_x M_0 \leq \mathbb E_x M_\tau =
\mathbb E_x u(x_\tau) + c \mathbb E_x \tau \epsilon^3 = \mathbb E_x
u(x_\tau) + O(\epsilon).$$ Thus, $u_1^\epsilon(x) \geq u(x) -
O(\epsilon)$.  A symmetric argument gives $u_2^\epsilon(x) \leq u(x)
+ O(\epsilon)$, so $|u_1 - u|_\infty = O(\epsilon)$, since
$u_1^\epsilon \leq u_2^\epsilon$.

\end{proof}

We now show that the solution $u$ of $\Delta_p u = g=-2f/\beta$ has a
natural game-theoretic interpretation as the limiting value of a
modified tug of war game in which in addition to payoff player I
receives when the game terminates, player I receives a ``running
payoff'' of size $f(x_k) \epsilon^2$ at the $k$th step of the game.

\begin{theorem} \label{runningpayofftheorem}
Suppose that $u$ is smooth on a bounded domain $\Omega'\supset
\overline \Omega$ with non-vanishing gradient and $\Delta_p u = g$,
where $f = -\frac{\beta}{2}g$ is bounded below on $\Omega$ by a
positive constant. Then $u$ is the limit as $\epsilon \to 0$ of the
functions $u^\epsilon_1$ defined for tug of war with noise on $\Omega$, with
running payoff $f$ and boundary payoff given by $F = u$ on $\partial
\Omega$.
\end{theorem}

\begin{proof}
The proof is exactly the same as in Theorem
\ref{pharmoniclimitsmoothcase} (the estimates in Lemma
\ref{onestepexpectationbound} include the running payoff case),
except that we need a different argument to show that the expected
number of turns in the game is $O(\epsilon^{-2})$.  But since $f$ is
bounded below, if player I adopts the strategy of always pulling
distance $\epsilon$ in the gradient direction, either the expected
number of turns in the game is less than
$$\frac{\sup_\Omega u(x) - \inf_{\partial \Omega} F(x)}{\inf f(x)}
\epsilon^{-2},$$ in which case the arguments of Theorem
\ref{pharmoniclimitsmoothcase} apply, or it is greater than this, in
which case player I's expected payoff is even greater than $u(x)$.

On the other hand, player II can adopt a strategy which makes the
process $u(x_k) + \sum_{i=0}^{k-1} f(x_i)\epsilon^2 - c k\epsilon^3$
a supermartingale (for some constant $c$ independent of $\epsilon$
and $x_0$).  Since this supermartingale is bounded below, by
optional stopping
$$\mathbb E \left[u(x_\tau) +\sum_{i=0}^{\tau-1} f(x_i)\epsilon^2 -
c \tau\epsilon^3 \right]\leq u(x_0).$$ Hence
$$\mathbb E \left[ u(x_\tau) + \tau \left((\inf_{x \in \Omega} f(x)\epsilon^2) - c \epsilon^3 \right)\right]\leq
u(x_0).$$ Since $u(x_\tau)$ is bounded on $\Omega$, this implies
$\mathbb E \tau = O(\epsilon^{-2})$.
\end{proof}

The remainder of this paper will focus on zero running payoff case
$f = 0$.

\subsection{Continuous boundary conditions, regular boundary}

Even if the boundary function $F$ is continuous and $u^{\epsilon}$
converges pointwise to a function $u$, it may not be the case that
$u$ extends continuously to $F$ on $\partial \Omega$.  For example,
in dimension two, if $p \leq 2$ and $\Omega$ is the unit ball minus
the origin and $f$ is the function which is $1$ at the origin and
$0$ at the boundary of the unit ball, then it is not hard to see
that as $\lim_{\epsilon \rightarrow 0} u^\epsilon_1 = 0$ throughout
$\Omega$, which is discontinuous at the origin.  However, we will
see that this cannot happen if $\Omega$ is sufficiently regular.
This section extends the results of the previous section to
sufficiently regular domains by proving Theorem
\ref{regularandcontinuousconvergence}.


\proofof{Theorem \ref{regularandcontinuousconvergence}}{\bf (i)} Suppose that
$\Omega$ is game-regular.  Our strategy will be to use this fact to
establish some a priori equicontinuity on the $u^\epsilon_1$ that
will ensure (via a compactness argument) subsequential uniform
convergence of the $u^\epsilon_1$ to a limit that extends
continuously to $F$ on $\partial \Omega$. Then we will use
comparison arguments and Theorem \ref{pharmoniclimitsmoothcase} to
show that any such limit is a viscosity solution to $\Delta_p u = 0$,
which is unique by Proposition \ref{pharmonicunique}.

Fix a continuous $F$ and write $b_L = \min_{y \in \partial \Omega}
F(y)$ and $b_U = \max_{y \in \partial \Omega} F(y)$. Given a
constant $\gamma$, we can find a $\delta$ such that $|F(x) - F(y)| <
\gamma/2$ whenever $|x-y| < \delta$.  We then define $\eta =
\frac{\gamma}{2(b_U- b_L)}$.

By game-regularity, for each $y \in
\partial \Omega$ we can find $\delta_0 = \delta_0(y)$ and $\epsilon_0 =
\epsilon_0(y)$ such that for every $x_0$ with $|x_0 - y| < \delta_0$
and $\epsilon < \epsilon_0$, player I has a strategy that guarantees
that an $\epsilon$-step game started at $x_0$ will terminate at a
point on $\partial \Omega \cap B(y,\delta)$ with probability at
least $1 - \eta$.  The probability that $F$ at the terminal point
differs from the $F(y)$ by more than $\gamma/2$ is at most
$\frac{\gamma}{2(b_U- b_L)}$, so we have that $|u^\epsilon_1(x) -
F(y)| \leq \gamma$ for each $x \in B(y,\delta_0)$.

By compactness of $\partial \Omega$ (a closed bounded set), we can
find a finite collection $y_1, \ldots, y_N$ such that $$\partial
\Omega \subset S = \cup_{i=1}^N B(y_i,\delta_0(y_i)/2).$$ Now,
define $$\overline \delta_0 = \text{dist}(\Omega \backslash S,
\partial \Omega) \wedge \min_{1\leq i \leq N} \delta_0(y_i)/2,$$
and $\overline{\epsilon}_0 = \min_{1 \leq i \leq N} \epsilon_0(y_i)
\wedge \overline \delta_0/(2\alpha)$. Every point $x$ within
distance $\overline \delta_0$ of $\partial \Omega$ has the property
that $|u^\epsilon_1(x) - u(x')| \leq 2\gamma$ for every $x'$ with
$|x-x'| \leq \overline \delta_0$ and $\epsilon < \overline
\epsilon_0$, because $x \in B(y_i,\delta_0(y_i)/2)$ for some $i$,
and hence $x' \in B(y_i,\delta_0(y_i))$. This in turn implies that
\begin{equation} \label{equi}
\mbox{ \rm for every $\epsilon < \overline \epsilon_0$ and any
$x, x' \in \Omega$ with $|x -
x'| < \overline \delta_0$, we
have } |u^\epsilon_1 (x) - u^\epsilon_1(x')| < 2 \gamma \, .
\end{equation}
 The reason
is that a player starting at $x'$ can always adopt the same strategy
that would have been chosen starting at $x$ (i.e., define $x'_j = x_j +
x'-x$ and when the game position is some $x'_{k-1}$ a player can
choose the step $v_k$ that would be chosen if the game position
were $x_{k-1}$) up until the first time that $\text{dist}(x_k,
\partial \Omega) \wedge \text{dist} (x_k', \partial \Omega)
< \overline \delta_0$.
At that point we have $|u^\epsilon_1(x_k') - u^\epsilon_1(x_k)| \leq
2 \gamma$, and the result follows. By symmetry, the analog of (\ref{equi})
with $u_1^\epsilon$ replaced by $u_2^\epsilon$ also holds.
Consider $\underline{u}_1(x):=\liminf_{\eps \to 0} u^\epsilon_1 (x)$
and $\overline u_2(x):=\limsup_{\eps \to 0} u^\epsilon_2 (x)$.
 The asymptotic equicontinuity (\ref{equi}) implies that
 $\underline{u}_1$ and  $\overline{u}_2$ are (uniformly)
 continuous functions on $\overline{\Omega}$ which agree on
$\partial \Omega$.
If we can show that these functions are both $p$-harmonic,
then uniqueness of extensions (Proposition \ref{pharmonicunique})
 and the general inequality $u_1^\epsilon \le u_2^\epsilon$
will imply that the pointwise limits $\lim_{\eps \to 0} u^\epsilon_i (x)$
for $i=1,2$ exist and coincide everywhere in $\Omega$.
(Another application of (\ref{equi}) will then yield that the
 convergence is uniform.)

Using symmetry, it only
 remains to show that $\underline{u}_1$ is $p$-harmonic.
By Proposition
\ref{viscosityweakvariationalequivalance} it is enough to show that
$\underline{u_1}$ is a viscosity solution of $\Delta_p u=0$ in $\Omega$.
We will verify that $\underline{u}_1$ is a viscosity subsolution;
the proof that it is a supersolution is similar.

Let $\phi$ be a $C^2$ function in a neighborhood $V_0$
of $x\in \Omega$ such
that $\phi(x) = \underline{u}_1(x)$, the inequality
$\underline{u}_1 \leq \phi$ holds in $V_0$ and $\nabla \phi(x) \not = 0$
in $V_0$. We must show that  $\Delta_p \phi(x) \geq 0$.
Otherwise, there exist $r>0$ and $\theta>0$ such that
 $\Delta_p \phi \le -\theta$ on $B(x,r) \subset V_0$.
Let
$$
\wphi(y)=\phi(y)+\frac{\theta}{4d} (y-x)^T(y-x)-\frac{\theta r^2}{8d} \,,
$$
denote $g=\Delta_p \wphi \le \Delta_p \phi+\theta/2 \le -\theta/2$,
and consider noisy tug of war played in $B(x,r)$ with stepsize
$\epsilon$, running payoff $f=-\frac{\beta}{2} g \ge \frac{\beta
\theta}{4}$ and boundary values $\wphi$ on $\partial B(x,r)$. By
Theorem \ref{runningpayofftheorem}, the value functions for player I
in this game (which must be greater than $u_1^\eps$ for small $\eps$
due to the positive running payoff and the larger boundary values),
 converge as $\epsilon \to 0$ to $\wphi$.
This is a contradiction, since $\wphi(x)< \overline{u}_1(x)$.
Thus we have shown that  $\Delta_p \phi(x) \geq 0$, and the proof is complete.

\medskip

\noindent{\bf Proof of
Theorem \ref{regularandcontinuousconvergence}(ii)}. \;
Suppose that for every continuous $F$, the
$u_1^\epsilon$ converge uniformly.  Then we claim that $y$ is game
regular, i.e., i.e., that for every $\delta> 0$ and $\eta > 0$ there
exists a $\delta_0$ and $\epsilon_0$ such that for every $x_0$ with
$|x_0 - y| < \delta_0$ and $\epsilon < \epsilon_0$, player I has a
strategy that guarantees that an $\epsilon$-step game started at
$x_0$ will terminate at a point on $\partial \Omega \cap
B(y,\delta)$ with probability at least $1 - \eta$.  To see this,
define $F(y') = -|y-y'|$ on $\partial \Omega$ and let $u$ be its
extension. Given $\delta$ and $\eta$ we can choose $\delta_0$ such
that $|x_0 - y| < \delta_0$ implies $|u(x_0)-u(y)| < \eta\delta/2$.
By uniform convergence, we can find $\epsilon_0$ such that for
$\epsilon < \epsilon_0$ we have $||u^\epsilon_1 - u||_\infty <
\eta\delta/2$ which implies that $|u_1^\epsilon(x_0)| <\eta\delta$,
and hence player one has a strategy that guarantees that the
probability of ending $\delta$ away from $y$ is at most $\eta$. \qed

\subsection{Sufficient conditions for regularity}
We first present a condition equivalent to game-regularity which is easier to verify.

\begin{lemma} \label{equivreg}
Fix  $p>1$, a measure $\mu$ as in the introduction and a
domain $\Omega \subset \R^d$. A point $y \in
\partial \Omega$ is  a game-regular boundary point of
$\Omega$ if and only if there is some $\theta>0$ with the following
property: for every $\delta > 0$ there exists a $\delta_0$ and
$\epsilon_0$, such that for every $x_0$ with $|x_0 - y| < \delta_0$
and $\epsilon < \epsilon_0$, player I has a strategy that guarantees
that an $\epsilon$-step game started at $x_0$ will terminate at a
point on $\partial \Omega$ before exiting $B(y,\delta)$ with
probability at least $\theta$.
\end{lemma}

\begin{proof}
Given $\theta$ as in the statement of the lemma and $\eta$ as in the definition,
find $k$ such that $(1-\theta)^k<\eta$.
For $\eps_0$ and $\delta_0$ as in the statement, there exist $\delta_1$ and $\eps_1$
such that for every
$x_0$ with $|x_0 - y| < \delta_1$ and $\epsilon < \epsilon_1$,
player I has a strategy that guarantees that an $\epsilon$-step game
started at $x_0$ will terminate at a point on $\partial \Omega $ before exiting
$B(y,\delta_0-\alpha \eps)$, with probability at least $\theta$.
Iterating this argument $k$ times proves the lemma.
\end{proof}

The following is needed for the proof of Proposition
\ref{coneimpliesgameregular}:

\begin{lemma} \label{annulusbound}
Fix $p$ with $1 < p \leq \infty$. Given positive constants $r,s,t,
\epsilon$ with $0<s<1<t$, consider an $\epsilon$ game of tug of war
with noise on an annulus $B(0,t r) \backslash \overline {B(0,s r)}$
with $1 < p \leq \infty$ and $0 < s < 1 < t$. If the initial game
position is some $x_0$ with $|x_0| - \alpha\epsilon \leq r \leq
|x_0|+\alpha\epsilon$, and the ratio $t/s$ is held fixed, then as
$\frac{\epsilon}{rs} \to 0$, player I has a strategy that guarantees
that the game position will terminate on $\partial B(0,sr)$ (instead
of $\partial B(0,tr)$) with probability at least $b -
O(\frac{\epsilon}{rs})$ where
$$b=b(s,t,p,d):= \frac{t^c - 1}{t^c - s^c}$$ where $c = c(p,d) :=
\frac{p-d}{p-1}$.  (As in Theorem \ref{pharmoniclimitsmoothcase},
the constant implied by $O(\cdot)$ may depend on $\mu$.)

\end{lemma}

\begin{proof}
Given $d \geq 1$ and $p>1$ and $x \in \mathbb R^d \backslash \{0\}$,
we write
$$\rho_{d,p}(x) := \begin{cases} |x|^{c(p,d)} & p \not = d \\
 \log |x| & p = d. \\ \end{cases},$$
where $c = c(p,d) = \frac{p-d}{p-1}$. The reader may check that for
each $1 < p < \infty$ the function $\rho_{d,p}(x)$ is a radially
symmetric $p$-harmonic function on $\mathbb R^d \backslash \{0\}$.
Taking $F$ to be the function that has these values on the boundary
of the annulus and applying Theorem \ref{pharmoniclimitsmoothcase},
 we see that if $r$, $s$, and $t$ are fixed,
player I can
achieve a probability of reaching the inside first that is, up to an
$O(\epsilon)$ error, given by
$$\frac{\rho_{d,p}(tr) - \rho_{d,p}(r)}{\rho_{d,p}(tr) -
\rho_{d,p}(sr)}= \frac{(tr)^c - r^c}{(tr)^c - (sr)^c} = \frac{t^c -
1}{t^c - s^c}.$$  (A similar argument applies when $p=d$.) If only
the ratio $t/s$ is fixed, then the fact that the error is
$O(\frac{\epsilon}{rs})$ follows similarly from Theorem
\ref{pharmoniclimitsmoothcase} and a rescaling that replaces
$\epsilon$ with $\epsilon' = \epsilon/(rs)$ and replaces $r$ with
$r' = r/(rs)$ (so that $r's = 1$ and $r't$ is fixed).

\end{proof}

\proofof{Proposition \ref{coneimpliesgameregular}}
 {\bf (i)} We first prove
that all domains are game-regular when $p > d$. Fix $s = 1/2$ and
$t=2$.  Then since $c(p,d) > 0$, we have $b = b(s,t,p,d) > 1/2$.
Consider the domain $\Omega = \R^d \backslash \{0 \}$.  For each
integer $m \geq 0$, write $r_m = 2^m \alpha \epsilon$ and let $A_m$
denote the closed annulus centered at $0$ with radii $r_m + \alpha
\epsilon$ and $r_m - \alpha \epsilon$.

By Lemma \ref{annulusbound}, if the game position begins at $x_0 \in
A_m$, and $m \geq 1$, then Player I can arrange, with probability at
least $b - O(\epsilon/r_m) = b - O(2^{-m})$ to have the game
position enter $A_{m-1}$ before entering $A_{m+1}$. Next, let $P_m$
denote the infimum over all starting points $x_0 \in A_m$ of the
maximum probability with which player I can guarantee that the game
position hits $\overline B(0, \alpha \epsilon)$ before hitting
$A_{m+1}$. It is not hard to see that $P_m > 0$ for every $m \geq
1$. (If player I wins a long sequence of coin tosses---and chooses a
length $\epsilon$ vector uniformly at random each time---then the
game position will undergo a symmetric random walk that moves by at
most $\alpha \epsilon$ at each step.  It is enough to observe that
any such random walk hits $\overline B(0, \alpha \epsilon)$ before
$A_{m+1}$ in some fixed number of steps with a probability that is
bounded below independently of the choice of starting position in
$A_m$.) Since player I can arrange to have at least a $b_m =
b-O(2^{-m})$ chance of reaching $A_{m-1}$ before $A_{m+1}$ when the
game positions starts in $A_m$, we have
$$P_m \geq b_m\Bigl( P_{m-1} + (1-P_{m-1}) P_m \Bigr).$$
Rearranging terms, we get
$$(1-b_m + b_mP_{m-1})P_m \geq b_m P_{m-1}$$
and hence
$$P_m \geq \frac{b_m P_{m-1}}{1 - b_m + b_m P_{m-1}} = \frac{P_{m-1}}{b_m^{-1}-1 + P_{m-1}}.$$
The right hand side is less than or greater than $P_{m-1}$ as
$P_{m-1}$ is, respectively, greater than or less than $2-b_m^{-1}$.
From this it is not hard to see that $P_m$ is bounded below
independently of $m$.  In fact, $\liminf P_m \geq 2-b^{-1}$.

Now, to prove game-regularity, we must show that there exists
$\theta > 0$ such that for every $\delta
> 0$, there exists a $\delta_0$ and $\epsilon_0$ with the property
that for every $x_0$ with $|x_0 - y| < \delta_0$ and any $\epsilon <
\epsilon_0$, player I has a strategy that guarantees that an
$\epsilon$-step game started at $x_0$ will terminate at a point on
$\partial \Omega \cap B(y,\delta)$ with probability at least
$\theta$. We may take $\delta_0 = \delta/4$, so that for every
sufficiently small $\epsilon$ we have $\delta_0 \leq r_m < r_{m+1}
\leq \delta$ for some $m$. Then if $|x_0 - y| < \delta_0$, player I
can arrange to reach $\overline B(y, \epsilon)$ (or terminate the
game sooner) before exiting $B(y,\delta)$ with probability at least
$P_m$, which is bounded below by some constant $\theta$ independent
of $m$.

\noindent{\bf (ii)} The second part of the proposition states that all domains having
the cone property are game-regular.  Assume without loss of
generality that $y = 0 \in \partial \Omega$ is the tip of the cone,
and we will argue that $y$ is game-regular. There exists some
constant $\chi \in (0,1)$ and $R$ such that for every $r<R$, there
exists a ball of radius $\chi r$ contained in $\mathcal C$ whose
center is of distance $r$ from the origin.  If the game position
begins at a point $x_0$ whose distance is $r$ from the tip of the
cone, then player I can adopt the strategy of pulling towards the
center of a ball of radius $\chi r$ centered at a point $z$ in
$\mathcal C$ (which is distance at most $2r$ away from $x_0$). Then
by Lemma \ref{annulusbound}, for all sufficiently small $\epsilon$
player one has a probability $\theta$ (bounded below independently
of $\epsilon$ and $r$) of reaching this ball and terminating the
game (or terminating the game even sooner) before exiting $B(0,4r)$.
The result now follows from Lemma \ref{equivreg}.

\noindent{\bf (iii)} Finally, we prove that all simply connected domains are game-regular
when $d=2$. By translation, we can suppose that $0 \in \partial \Omega$.
Our goal is to find $\theta>0$, so that for any initial game position $x_0$, for
sufficiently small $\epsilon$ player I can ensure that the game will
end with probability at least $\theta$
before the game position reaches a point outside of $B(0,2|x_0|)$.

By rotating and scaling, we may and shall assume that
$x_0 = (-2,3)$.  Let $L_1$ be the ordered line segment in $\R^2$ from
$(-2,3)$ to $(1,-1)$.  Let $L_2$ be the ordered line segment from
$(1,-1)$ to $(-1,-1)$ and let $L_3$ be the ordered line segment from
$(-1,-1)$ to $(2,3)$.  Let $L$ be the concatenation of $L_1$, $L_2$,
and $L_3$ (a continuous path). Finally, let $z_0 = (-1,2), z_1, z_2,
\ldots, z_{1200} = (1,2)$ be an evenly spaced sequence of points along
the path $L$, and let $r = \frac{|L_1| + |L_2| + |L_3|}{1200} =
\frac{12}{1200} = 1/100$ be the distance between adjacent $z_i$
along the path $L$.  If the game position begins at a point in
$B(z_i,r) \subset B(z_{i+1},2r)$ for $0 \leq i < 1199$, then for
any $b'$ with $0 < b' < b(1/2,2,p,d)+O(r^{-1})$ player I can arrange
(for all $\epsilon$ sufficiently small) to reach a point in
$B(z_{i+1},r)$ before exiting $B(0,4r)$ with probability at least
$b'$. Thus, if the game position begins at $z_0$, player I can
arrange, with probability $\theta = (b')^{1200}$, to hit each of the
$B(z_i,r)$ without ever reaching a point of distance more than $4r$
from $L$.
\begin{figure}[htbp]\label{fig2}
\begin{center}
\epsfig{file=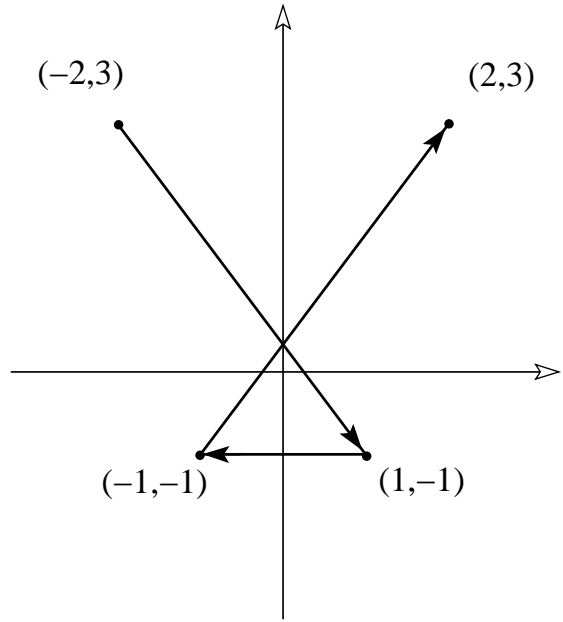, width=.45 \textwidth}
\caption{ Path used to prove simply connected domains are game-regular}
\end{center}
\end{figure}
Since $4r<\dist(L,0)$, a simple topological argument (this is where
$d=2$ is used) shows that if $0$ lies on the boundary of $\Omega$
and $\Omega$ is simply connected, then the path obtained by joining
subsequent game positions $x_j$ in such a game sequence with line
segments would have to surround $0$ and intersect itself---and would
thus have to cross $\partial \Omega$ since $\Omega$ is simply
connected. However, this also implies that at some point the game
position is within $\alpha \epsilon$ of $\partial \Omega$, and thus
the game must have terminated. We conclude that  for all
sufficiently small $\epsilon$, player I can arrange for the game to
end with probability at least $\theta = (b')^{1200}$ before the game
position reaches a point outside of $B(0,2|x_0|)$. The result
follows from Lemma \ref{equivreg}.

\qed

\section{Variants} \label{variantsection}
\subsection{Alternating turns}

We now introduce a two-player zero-sum game, called {\bf alternating
turn tug of war (with noise)}.  The rules are exactly the same as
those of ordinary tug of war except that players alternate turns
(with player I moving first) instead of determining turn order with
coin tosses.

In the alternating turn game, we define $p = p(\mu) \in [1, \infty]$
by $p = \frac{C_{1,1}+ C_{2,2}}{C_{2,2}}$ instead of $p =
\frac{C_{1,1}+C_{2,2}+1}{C_{2,2}}$.

Equivalently, $p$ is such that for some $\beta
> 0$, we have $C_{1,1} = \beta q^{-1}$ and $C_{i,i} = \beta p^{-1}$
for $i \geq 2$ for some $\beta > 0$, where $p^{-1} + q^{-1} = 1$. As
in the random turn game, $p$ is chosen so that if one player always
chooses the vector $v_k$ to be some $v$ with $|v|=\epsilon$, and the
other player always chooses $v_k$ to be $-v$, then for each $k$ the
variance of $x_k - x_{k-1}$ is proportional to $p\epsilon^2$ in the
direction of $v$ and $q\epsilon^2$ in each direction orthogonal to
$v$.

As in the random turn case, we let $\phi(x) = x^TAx + (\xi,x)$ (for
$\xi \in \R^d$) be a quadratic function. Suppose that $x_0 = 0$ and
let $\psi(v)$ be the expected value of $\phi(x_1) = \phi(v_1 + z_1)$
if player I chooses $v_1 = v$.  Then $z_1$ has expectation zero and
covariance matrix $C$ given by
$$\beta q^{-1} v v^T + \beta p^{-1} (|v|^2I - vv^T).$$

Since $z_1$ has expectation zero, we have $$\psi(v) := \mathbb E
[\phi(v+z_1)] = \phi(v) + \sum_{i,j} A_{i,j} \mathbb E [z_1^i
z_1^j],$$ where $z_1^i$ is the $i$th component of $z_1$).

We can also write this as $\sum_{i,j} A_{i,j} C_{i,j}$. Note that
$\sum (v v^T)_{i,j} A_{i,j} = v^T A v$ and $\sum I_{i,j} A_{i,j} =
\text{Tr} A$.  We can now write

\begin{eqnarray*}
\psi(v) &= &\phi(v) + (\beta q^{-1}) v^T A v + \beta p^{-1}
(|v|^2\text{Tr}
A - v^T A v) \\
&= &\phi(v) + (\beta q^{-1} - \beta p^{-1})v^TAv + \beta p^{-1} |v|^2\text{Tr}A\\
& =& (\beta q^{-1} - \beta p^{-1} + 1) v^TAv + \beta p^{-1}
|v|^2\text{Tr} A + (\xi, v)
\end{eqnarray*} and we can write $\psi(v) = (\xi, v) + v^T B v$ where $B =
(\beta q^{-1} - \beta p^{-1} + 1)A + \beta p^{-1}(\text{Tr}A) I$.

We give the following alternating turn analog of Lemma
\ref{onestepexpectationbound}.

\begin{lemma} If $\phi$ is quadratic, then for any even $k$ (in the alternating turn game)
if player I makes the $k+1$th move in the gradient direction, then
regardless of what player II's strategy is, we have
$$\mathbb E [\phi(x_{k+2})|h_k] \geq \phi(x_k) + \beta \Delta_p
\phi(x_k)\epsilon^2 + O(\epsilon^3).$$ The constant in the
$O(\epsilon^3)$ depends only $\|B\|$ and $|\xi| = |\nabla
\phi(x_k)|$, as in Lemma \ref{maximalgainperturn}.
\end{lemma}

\begin{proof}
It is enough to prove this for $k=0$.  We may assume $x_0=0$ and
then with the gradient strategy $x_1 = \epsilon|\xi|^{-1}\xi$ and
$\nabla \phi(x_1) = \xi + 2Ax_1$.  Therefore $$|\nabla\phi(x_1)|^2 -
|\nabla \phi(x_0)|^2=4\eta^TAx_1 + O(\epsilon^2).$$  Dividing by
$|\nabla \phi(x_0)| + |\nabla \phi(x_1)| = 2|\xi| + O(\epsilon)$
gives $$|\nabla \phi(x_1)| - |\nabla \phi(x_0)| = \frac{2 \xi^T A
x_1}{|\xi|} +O(\epsilon^2) = \epsilon \Delta_\infty \phi(x_0) +
O(\epsilon^2).$$

If player I is moving at the first turn, player I can arrange---by
pulling directly in the gradient direction---to have

\begin{equation}\label{onestepestimate} \mathbb E [\phi(x_1)] - \phi(x_0) \geq |\nabla
\phi(x_0)|\epsilon + \frac{\epsilon^2}{2} \left( \beta\Delta_p +
\Delta_\infty \right) \phi(x_0) + O(\epsilon^3).\end{equation} When
player II makes the subsequent move, the expectation is smallest (up
to $O(\epsilon^3)$ error) if player II always pulls in the minus
gradient direction. In that case we have

\begin{eqnarray*}\mathbb E [\phi(x_2)|h_1] - \phi(x_1) &\geq& -|\nabla
\phi(x_1)|\epsilon + \frac{\epsilon^2}{2} \left( \beta\Delta_p +
\Delta_\infty \right) \phi(x_1) + O(\epsilon^3) \\
&=& -\Bigl(|\nabla \phi(x_0)|  +\epsilon \Delta_\infty \phi(x_0)
\Bigr)\epsilon + \frac{\epsilon^2}{2} \left(\beta\Delta_p +
\Delta_\infty \right) \phi(x_0) + O(\epsilon^3) \, .\end{eqnarray*}

Taking expectations and adding (\ref{onestepestimate}) gives

$$\mathbb E \phi(x_{2}) - \phi(x_0) \geq \beta \Delta_p
\phi(x_0)\epsilon^2 + O(\epsilon^3) \, ,$$ as claimed.
\end{proof}

The proofs of analogs of Theorem \ref{pharmoniclimitsmoothcase},
Theorem \ref{runningpayofftheorem}, and Theorem
\ref{regularandcontinuousconvergence} are exactly the same as in the
random turn case.

\subsection{Direction selection and Spencer's game}

The following game (along with other variants) was introduced by
Spencer \cite{MR0526057} and studied in detail by Kohn and Serfaty
\cite{MR2200259}.  We focus on the case $d=2$, although many natural
variants exist in higher dimensions.

Fix $\Omega, x, F$ as in tug of war and set $x_0 = x$. At each turn
$k$, player II chooses a vector $v_k$ of length exactly $\epsilon$.
Player I then chooses a sign $\sig_k \in \{-1, 1\}$ and sets $x_k =
x_{k-1} + \sig_k v_k$.  Kohn and Serfaty were primarily concerned
with convex domains, constant running payoff $f$ and zero boundary
conditions $F=0$ (i.e., one player seeks to maximize the number of
steps before the boundary of the domain is reached and the other
player seeks to minimize that number).  In this case, the relevant
operator is $\Delta_1$, i.e., we generally expect the limiting
payoff function $u$ to be a solution to $\Delta_1 u = -2f$; however,
many complicating issues arise in this setting that do not arise in
tug of war, so the results here are somewhat more restrictive. Even
on simple domains (such as the unit disc) with $C^\infty$ boundary
values $F$, the equation $\Delta_1 u = 0$ may not have a unique
solution.  (For a concrete example, note that if the unit circle is
parameterized by angle $\theta \in [0, 2\pi)$, then the function
$\cos(2\theta)$ can be extended to the interior of the unit disc in
such a way that it is constant on vertical lines; it can also be
extended to the unit disc in such a way that it is constant on
horizontal lines. Both extensions are solutions to $\Delta_1 u = 0$
--- at least in the sense that they are smooth and satisfy $\Delta_1 u = 0$ wherever
$\nabla u \not = 0$.)  We refer the reader to \cite{MR2200259} for
more discussion along these lines.

Consider the following interpolation between tug of war and
Spencer's game:  as before, we set an initial game position $x_0 \in
\Omega$.  At the $k$th turn, we toss a $p^{-1}$-coin to determine
how the game position is updated: with probability $p^{-1}$, $x_k$
is updated using the rule of Spencer's game described above (with
player II choosing a vector and player I choosing the sign); with
probability $q^{-1}$, it is updated using the rule of ordinary tug
of war (i.e., a fair coin is tossed, and the player who wins the
toss chooses a $v_k \in \overline B(0,\epsilon)$ and we set $x_k =
x_{k-1} + v_k$).

In this setting, the analog of Theorem \ref{runningpayofftheorem} is
even simpler to prove than in the other two cases, since we may deal
with the two cases separately.  First, if it is decided that the
game position will be updated at the $k+1$th step using the tug of
war rule, then by pulling in the gradient direction, player I can
arrange for $\mathbb E[\phi(x_{k+1}) - \phi(x_k) \, | \, h_k]$ to be at least
$\Delta_\infty \phi/2 + O(\epsilon^3)$ when $\phi$ is smooth (and
player II can arrange for $\mathbb E[\phi(x_{k+1}) - \phi(x_k)\, | \, h_k]$ to
be at most $\Delta_\infty \phi/2 + O(\epsilon^3)$).

This follows as a special ($\mu$ supported at the origin) case of
the bounds given for tug of war with noise. Second, we claim if it
is decided that the game position will be updated at the $k+1$th
step using Spencer's rule, it is similarly clear (using the Taylor
expansion for $\phi$) that player II can arrange (by choosing the
direction orthogonal to the gradient) for $\mathbb E[\phi(x_{k+1}) -
\phi(x_k)\, | \, h_k]$ to be at most $\Delta_1 \phi/2 + O(\epsilon^3)$. If
player I adopts the strategy of always choosing a sign so that the
inner product $(\sig_k v_k, \nabla \phi(x_k))$ is non-negative, then
a similar Taylor expansion shows that $\mathbb E[\phi(x_{k-1}) -
\phi(x_k)\, | \, h_k]$ will be at least $\Delta_1 \phi/2 + O(\epsilon^3)$.
Thus, before it is decided which update rule will be used, either
player can arrange for the expected change, given the history, to be $\Delta_p \phi/2 +
O(\epsilon^3)$.  The remainder of the argument is the same as the
one given for ordinary tug of war.

A complete proof of an analog of Theorem
\ref{regularandcontinuousconvergence} in this context is also
possible, but the definition of game
regularity must be symmetrized to require the same for player II
as is required for player I.

\section{Irregular domains and discontinuous payoff functions}
\label{discontinuousboundarysection}

\subsection{Resolutive functions}

If either $\Omega$ is game irregular or $F$ is discontinuous (e.g.,
if $\Omega$ is the unit disc and $F$ is zero on points with rational
angles and one on all other points), then we cannot expect there to
exist a $p$-harmonic function on $\Omega$ that extends continuously
to $F$ on $\partial \Omega$. There is still a natural notion of
$p$-harmonic extension of $F$ provided $F$ is {\em resolutive} (as
defined below; when $p=2$, resolutivity follows from boundedness and
Borel measurability). However, in this setting, the $u^\epsilon$
defined above need not converge to that extension (as is clear in
the example where $\Omega$ is the unit disc and $F$ is the
characteristic function of the set of points of irrational angles).

To generalize $p$-harmonic extensions to general domains and
boundary conditions, it is conventional to define the {\bf upper
class} $\mathcal U_F$ of $F$ as the set of all functions $u$ such
that $u$ is $p$-superharmonic in $\Omega$, $u$ is bounded below, and
$\lim \inf_{x \rightarrow y} u(x) \geq F(y)$ for all $y \in
\partial \Omega$.  The {\bf lower class} $\mathcal L_F$ can be
analogously defined by writing $v \in \mathcal L_F$ if $-v \in
\mathcal U_{-F}$.  The function $\overline H_F = \inf \{u : u \in
\mathcal U_F \}$ is the {\bf upper Perron solution} of $F$ in
$\Omega$ and $\underline H_F = \sup \{u: u \in \mathcal L_F \}$ is
the {\bf lower Perron solution} of $F$ in $\Omega$.  The upper and
lower Perron solutions are either identically $\pm \infty$ or
everywhere finite and $p$-harmonic. In particular, when $F$ is
bounded, both $\underline H_F$ and $\overline H_F$ are bounded and
$p$-harmonic (Chapter 9 of \cite{MR1207810}).

We say that $F$ is {\bf resolutive} if $\underline H_F$ and
$\overline H_F$ agree and are $p$-harmonic in $\Omega$.  If $\Omega$
is regular then every continuous function on $\Omega$ is resolutive;
in fact, if $\Omega$ is regular, then every bounded and lower
semi-continuous $F$ on $\partial \Omega$ is resolutive (Chapter 9 of
\cite{MR1207810}). If $F_i$ are resolutive and $F_i \rightarrow F$
uniformly then $F$ is resolutive.  It is known that when $p=2$ any
bounded Borel measurable $f:\partial \Omega \rightarrow \mathbb R$
is resolutive.  For general $p$, no such comparable result is known,
and it is an open problem whether every Borel function is resolutive
\cite{MR1207810}.  It is not hard to see that some simple
discontinuous functions (e.g., the function that is $1$ on points of
the boundary with rational coordinates and $0$ on all other points,
etc.) are resolutive for all $p$.

\subsection{Shrinking step sizes}
Consider the variant of tug of war called {\bf small-step-size tug
of war} played as follows: fix parameters $\Omega$, $F$, $p$, and
$x$ as in ordinary (either random turn or alternating turn) tug of
war. One player chooses any real $\epsilon_0>0$ and the other player
chooses any real $\epsilon$ with $0 < \epsilon \leq \epsilon_0$. The
players then play tug of war using using the parameter $\epsilon$.
Clearly, the payoff of this game is $\liminf_{\epsilon \rightarrow
0} u^\epsilon$ if player I chooses $\epsilon_0$ and player II
chooses $\epsilon$ and $\limsup_{\epsilon \rightarrow 0} u^\epsilon$
if player II chooses $\epsilon_0$ and player I chooses $\epsilon$.
The convergence in Theorem \ref{regularandcontinuousconvergence}
implies the following:

\begin{cor} \label{smallepsiloncorollary}
If $\Omega$ is game-regular, then for every continuous function $F:
\partial \Omega \rightarrow \mathbb R$, the value of
small-step-size tug of war as a function of the initial game
position is given by the unique $p$-harmonic function $u$ that
extends continuously to $F$ on $\partial \Omega$.
\end{cor}

Let $\underline v$ be the value for player I of the following game
called {\bf shrinking-step-size tug of war}.  Players begin by
choosing an $\epsilon_0$ and $\epsilon$ as in small-$\epsilon$ tug
of war, with the requirement that $d(x_0, \partial \Omega) >
\alpha\epsilon$ (so that the game does not end instantly).
Subsequently, at each time $j$ such that either
\begin{enumerate}
\item $d(x_j, \partial \Omega) \leq \alpha\epsilon$, or
\item for some $m\geq 1$, $j$ is the smallest integer for which
$d(x_j, \partial \Omega) \leq 2^{-m}$,
\end{enumerate} player I chooses a new $\epsilon_0$ and player II
chooses a new $\epsilon < \epsilon_0$ (again, with the requirement
that $d(x_j, \partial \Omega) > \alpha\epsilon$), and the game
continues. Clearly this game cannot terminate in finitely many
steps, because the $\epsilon$ used during a game step is always such
that the distance from the boundary is at least $\alpha \epsilon$
(and hence it is not possible to reach the boundary and end the game
during that step). We define the payoff to be the infimum of $F(y)$
over all limit points $y \in
\partial \Omega$ of the $x_j$, and $0$ if no limit point on $\partial \Omega$ exists.
Let $\overline v$ be the value for player II of the game defined
analogously except that player I chooses the $\epsilon_0$ and player
II chooses the $\epsilon$ (at each stage) and the payoff is the
supremum of $F(y)$ over all limit points of $y \in
\partial \Omega$, instead of the infimum. Clearly $\underline v
\leq \overline v$. Now we have the following:

\begin{theorem} \label{shrinkingstep}
In shrinking-step-size tug-of-war on a domain $\Omega$ with bounded
boundary function $F$, as defined above, the following holds
throughout $\Omega$:
$$\underline H_F \leq \underline v \leq \overline v \leq \overline
H_F.$$ In particular, if $F$ is bounded and resolutive, then its
$p$-harmonic extension $\underline H_F = \overline H_F$ is the value of
the game independently of which player chooses the $\epsilon_0$ and
whether the payoff is taken to be the supremum or the infimum of $F$
on the set of limit points the $x_k$.
\end{theorem}

We will make some more observations before proving this result.
When $F$ is the characteristic function $1_A$ for some $A \subset
\partial \Omega$, the value $\overline H_F(x)$ is called the {\bf
$p$-harmonic measure of $A$ at $x$} and written
$\omega_p(A,x,\Omega)$. It is well known (see, e.g.,
\cite{MR1207810}) that if $\omega_p(A,x_0,\Omega)=0$ for some $x_0$
in a connected domain $\Omega$, then $\omega_p(A,x,\Omega)=0$ for
all points $x\in \Omega$.

Theorem \ref{shrinkingstep} is interesting in light of the fact that
$p$-harmonic measure is non-additive even on null sets
\cite{MR2163560}. In fact, \cite{MR2163560} exhibits a disjoint
finite collection $\{A_i\}$ of resolutive sets with $p$-harmonic
measure zero whose union is all of $\partial \Omega$.

Interpreted game theoretically, this means that there can be a
partition $\{A_i\}$ of $\partial \Omega$ such that for each $i$,
there is a strategy that causes all of the limit points of the game
play to lie outside of $A_i$ with probability one.  In other words,
player I has a strategy for avoiding any {\em one} of the $A_i$ with
probability one (regardless of player II's actions), even though
player I clearly has no strategy for avoiding {\em all} of the $A_i$
with even positive probability.

Following Chapter 2 of \cite{MR1207810}, given a compact subset $K$
of $\Omega$, let $W(K, \Omega) := \{u \in C_0(\Omega): u \geq 1
\text{ on } K\}$ and define its {\bf $p$-capacity} to be
$\text{cap}_p(K, \Omega) = \inf_{u \in W(K, \Omega)} \int_\Omega
|\nabla u|^p$.  If $U$ is open, define $\text{cap}_p(U, \Omega)$ to
be supremum of $\text{cap}_p(K, \Omega)$ over compact $K \subset U$;
if $E$ is arbitrary, then $\text{cap}_p(E, \Omega)$ is the infimum
of $\text{cap}_p(U, \Omega)$ over open sets $U \supset E$.  From
Chapter 9 of \cite{MR1207810} we cite the following:

\begin{prop} The set of $p$-irregular boundary points has $p$-capacity
zero.
\end{prop}

Note also that a subset $A$ of $\Omega$ has {\bf $p$-capacity} equal
to zero if and only if $A \subset
\partial (\Omega \backslash A)$ (i.e., $A$ has empty interior) and
$\omega_p(A,x,\Omega)=0$ for each $x \in \Omega \backslash A$
\cite{MR1207810}. By Theorem \ref{shrinkingstep},
$\omega_p(A,x,\Omega)$ is an upper bound on the value function for a
game with payoff $1$ on $A$ and $\partial \Omega \setminus A$---with
equality in the case that $A$ is resolutive (which holds, in
particular, if $A$ is open or closed; see Chapter 9 of
\cite{MR1207810}).

Thus, Theorem \ref{shrinkingstep} gives a new interpretation of what
positive $p$-capacity means: a resolutive set has positive
$p$-capacity if and only if for some $x \in \Omega \backslash A$,
this value is non-zero.  In other words, the resolutive subsets of
$\Omega$ with positive $p$-capacity are precisely those sets $A$
that a player can reach---with positive probability---from some $x
\in \Omega \backslash A$ in shrinking-step-size tug of war.

One definition of a {\bf $p$-regular point} is a point such that
$\lim_{x \rightarrow x_0} \overline H_F(x) = F(x_0)$ for each
continuous $F:
\partial \Omega \rightarrow \mathbb R$.  A point is {\bf irregular} if it
is not regular.  It easy to see that game-regularity implies
$p$-regularity (since Theorem \ref{regularandcontinuousconvergence}
defines a $p$-harmonic extension explicitly in this case), but the
converse is not known. It may be that equivalence of game-regularity
and $p$-regularity depends very sensitively on the precise
termination rule used near the boundary. For this reason, we
consider the shrinking-step game to be more natural than the
ordinary game when the domain is wildly irregular.

\proofof{Theorem \ref{shrinkingstep}} Fix a bounded Borel measurable
function $F$ on the boundary of an arbitrary bounded domain $\Omega
\subset \mathbb R^d$.  Suppose that $u$ belongs to $\mathcal U_F$,
i.e., the set of all functions $u$ such that $u$ is
$p$-superharmonic in $\Omega$, $u$ is bounded below, and $\lim
\inf_{x \rightarrow y} u(x) \geq F(y)$ for all $y \in \partial
\Omega$.

It is enough to show that $u$ is an upper bound bound on the value
of shrinking-step-size tug-of-war for player one.  Let $\Omega_m$ be
the component of $\{x \in \Omega: d(x, \partial \Omega) > 2^{-m}\}$
containing $x_0$ (which is non-empty for all sufficiently large $m$)
and let $\tilde \Omega_m$ be any smooth domain with $\Omega_m
\subset \tilde \Omega_m \subset \Omega_{m+1}$. Let $u_m$ be the
$p$-harmonic extension to $\tilde \Omega_m$ of the values of $u$ on
$\partial \tilde \Omega_m$.  Then we know (recall Propositions
\ref{pharmonicunique} and \ref{viscosityweakvariationalequivalance})
that this extension is unique and that $u_m \leq u$ on $\Omega_m$.
By Theorem \ref{pharmoniclimitsmoothcase}, Player II can choose
$\epsilon$ small enough to ensure that the expected value of
$u(x_{k_m})$, where $k_m$ is the first $k$ for which which $x_k \not
\in \tilde \Omega_m$, is at most $u(x)$ (up to an error that can be
made arbitrarily small, say smaller than $\delta 2^{-m}$ for some
fixed $\delta > 0$).  This implies that $u(x_{k_m}) - \delta
\sum_{i=1}^m 2^{-m}$ is a supermartingale in $m$; hence $u(x_{k_m})$
almost surely converges to a limit $V$ with $\mathbb E V \leq u(x_0)
+ \delta$.  Since the payoff for player I is at most $V$, this
implies that the expected payoff for player I is at most $u(x_0) +
\delta$. \qed

\begin{cor}
If $F$ is resolutive then the value of $F$ on any subset of
$\partial \Omega$ of $p$-capacity zero of $\partial \Omega$ makes no
difference to the value of shrinking-step-size tug of war.  Zero
capacity sets cannot be reached with positive probability---i.e., if
$F$ is $1$ on a set of $p$-capacity zero and $0$ on the rest of
$\partial \Omega$, then $u=0$ identically.
\end{cor}

\noindent{\bf Remark.} As noted by the referee, similar results were proved by
analytical means in the papers by Kurki~\cite{kurki} and Avil\'es-Manfredi~\cite{AM}.

\subsection{Harmonic measure for $\Delta_p$ and porous sets}
\label{harmonicmeasuresection}

Here, we present an estimate for the $p$-harmonic measure of a porous set for $1 < p
< \infty$, analogous to the estimate given for $\omega_\infty$ in
\cite{pssw2}.  Estimates of this type were derived analytically by O. Martio,
see Theorem 11.27 in \cite{MR1207810}. Neverthless, we feel that the argument below shows
how useful game theoretic intuition is for these problems.

Recall that a set $S$ in a metric space $Z$ is
$\lambda$-\textbf{porous} if for every $r\in(0,\diam Z)$ every ball
of radius $r$ contains a ball of radius $\lambda\,r$ that is
disjoint from $S$. An example of a porous set is the ternary Cantor
set in $[0,1]$.

\begin{theorem} \label{poroustheorem}
Let $\Omega \subset \R^d$, $d>1$ be the closed unit ball. Let
$\lambda\in (0,1/2)$ and let $\delta>0$. Let $S$ be a
$\lambda$-porous subset of $\partial \Omega$, and let $S_\delta$ be
the closure of the $\delta$-neighborhood of $S$. Then for some
constant $c = c(\lambda, p) > 0$ we have
$$
\omega_p(S_\delta,0,\Omega) \le \delta^{c}.
$$
\end{theorem}

The key to the proof of the analogous result in \cite{pssw2} is that
in tug of war without noise, if  player I chooses a strategy of always
pulling towards a point $z$ at distance $r$ from the initial game
position $x_0$, and $z \in B(z,R) \subset \Omega$, then the limit as $\epsilon \to 0$ of the probability that
the  game sequence $\{x_j\}$  will reach  $z$
before  exiting $B(z,R)$ is  $(R-r)/R$; in particular, this limiting probability depends only on the
ratio $r/R$.

In our setting, this does not hold when $p < d$.  However, if
instead of targeting $z$ one merely aims to reach the ball
$B(z,r_0)$ for some  $r_0 < r < R$, then the (small $\epsilon$)
limit of the probability of reaching this ball is a positive
constant depending only on $r/r_0$ and $R/r_0$.

The idea of the proof in \cite{pssw2} was to construct a finite
sequence $x_0, \ldots, x_k$ of points---with $x_0=z$ and $x_k$ a
boundary point in the complement of $S_\delta$---such that the
distance between successive $x_i$ most some constant $r$ and the
distance from any $x_i$ to $S$ is at most some larger constant $R$.
Then player I can, with positive probability (call it $C$), reach
$x_1$ before exiting $B_R(x_1)$. Subsequently, player I can, with
probability $C$, reach $x_2$ before exiting $B_R(x_2)$.  Repeating,
we see that with probability at least $1-(1-C)^k$, player I can win
the game before exiting $\cup B_R(x_i)$.  The proof in \cite{pssw2}
then shows that if player one fails at this, then the same strategy
may be repeated on a smaller scale.  The proof of Theorem
\ref{poroustheorem} is identical to the one given in \cite{pssw2}
except that instead of targeting a sequence of individual points
(terminating at a point in one of the intervals of the complement of
$S_\delta$) one targets a sequence of small-radius balls
(terminating at a ball that lies just outside $\Omega$ and is
incident to the center of one of the intervals of the complement of
$S_\delta$).

\section{Non-measurable strategies and finite additivity} \label{valueexistencesection}

Some of the measurability restrictions we have imposed on strategies
and on $F$ can be relaxed by invoking some less standard mathematics
(such as finitely additive probability and/or the Axiom of
Determinacy).  There is a sizable stochastic game theory literature
on these subjects (see, e.g., \cite{MR1665779, MR1645524} and the
references therein) and we will only briefly mention what the main
results imply in our setting.

First, the requirement that the players adopt measurable strategies
is unnecessary if we require that all subsets of $\mathbb R^d$ be
$\mu$ measurable (since the latter implies that all strategies are
measurable).  This can be achieved within Zermelo Fraenkel set
theory with the Axiom of Choice (ZFC) by relaxing the requirement
that $\mu$ be $\sigma$-additive.  In this context, \cite{MR1665779,
MR1645524} define the expected payoff for {\em every} pair of
strategies (mixed or pure), and the following is immediate from
their more general result, which holds in ZFC:

\begin{prop} \label{generalvalueexistenceproposition}
Let $\Omega$ be an arbitrary domain, $x \in \Omega$, and $F$ a Borel
measurable payoff function.  Suppose we require that the $\mu$ in
the definition of tug-of-war is finitely additive and defined on
{\em all} subsets of $\mathbb R^d$, and that the players are allowed
arbitrary finitely additive mixed strategies. Then all of the
variants of tug of war we have introduced (alternating turn and
random turn; fixed $\epsilon$ and small-step-size and
shrinking-step-size) have values $u^\epsilon$.  Moreover,
$u^\epsilon_1 \leq u^\epsilon \leq u^\epsilon_2$, where
$u^\epsilon_1$ and $u^\epsilon_2$ are defined (as in the
introduction) using only pure measurable strategies.
\end{prop}

Recall that even for continuous boundary data, we did not resolve
the issue of when $u^\epsilon_1 - u^\epsilon_2$ under the
restriction to measurable strategies (although we showed that $|u^\epsilon_1 = u^\epsilon_2|_\infty \to 0$
as $\epsilon \to 0$).  Proposition
\ref{generalvalueexistenceproposition} gives a value for {\em every}
finite $\epsilon$ and every Borel boundary function $F$, which is at
least compatible with our measurable-setting definition in that it
is bounded between $u^\epsilon_1$ and $u^\epsilon_2$.

The results of \cite{MR1665779, MR1645524} also imply the existence
of a value function for shrinking step size tug of war, and it is
not hard to see (using comparison with smooth functions) that any
such value function must be $p$-harmonic.  We may interpret this
function as a canonical $p$-harmonic extension of $F$---defined for
every Borel measurable function $F$---(although we have not shown
that it does not depend on $\mu$).

An even more exotic assumption---which allows us to use arbitrary
payoff functions $F$---is the Axiom of Determinacy, which
contradicts the Axiom of Choice but which implies that all
stochastic games (in particular, tug of war games) have values for
all bounded payoff functions on $H_\infty$ \cite{MR1665779}. (It
also implies that that all subsets of $\mathbb R$ are Lebesgue
measurable \cite{MR1665779}.) Under the Axiom of Determinacy, the
value of shrinking-step-size tug of war exists (and is $p$-harmonic)
for {\em every} bounded boundary function $F:
\partial \Omega \rightarrow \mathbb R$.

\section{Open questions}

We conclude with some simple open questions.  We learned of the
following question from \cite{MR1207810}:

\begin{question} Are all Borel measurable functions $F$ of $\partial
\Omega$ resolutive when $\Omega$ is a bounded domain in $\mathbb
R^d$?
\end{question}

If the answer is no, then one may still attempt to construct a
canonical $p$-harmonic extension of $F$ by affirmatively answering
the following:

\begin{question} Does shrinking step size tug of war have a
value (within standard ZFC, using measurable strategies) for all
Borel measurable $F$ that depends only on $p$ (and not on any other
properties of $\mu$)?
\end{question}

We showed that game-regularity implied $p$-regularity but could not
prove the converse.

\begin{question} Is game-regularity equivalent to $p$-regularity?
\end{question}

In \cite{pssw2}, a game theoretic argument based on tug of war was
used to prove uniqueness of solutions to $\Delta_\infty u = g$
(given zero boundary conditions) for all sufficiently regular and
strictly positive $g$---and to show that uniqueness failed in
general if $g$ was allowed to assume values of both signs. It is
natural to ask whether these arguments can be adapted to the
$p$-Laplacian.

\begin{question} Let $\Omega$ be a bounded, game-regular domain and
$g$ a Lipschitz function on $\Omega$ with $\inf_{x \in \Omega} g(x)
> 0$.  Does there necessarily exist a unique solution to to $\Delta_p u =
g$?  Does this uniqueness fail when $g$ is allowed to assume values
of both signs?
\end{question}

Theorem \ref{regularandcontinuousconvergence} shows that
$u_1^\epsilon$ and $u_2^\epsilon$ converge to the same limit, but we
have not addressed the following:

\begin{question}
In the setting of Theorem \ref{regularandcontinuousconvergence}, is
it always the case that $u_1^\epsilon = u_2^\epsilon$?
\end{question}

\noindent{\bf Acknowledgement}.
We are grateful to Oded Schramm and David Wilson for sharing their insights on tug of war with us, and to 
Craig Evans for helpful discussions on the $p$-Laplacian and viscosity solutions. We thank
the referees for a careful reading of the manuscript and for several corrections and improvements.

\bibliographystyle{abbrv}
\bibliography{noisy}

\end{document}